\newfont{\Bbb}{msbm10 scaled\magstephalf}
\documentclass[reqno]{amsart}
\usepackage{amsfonts}
\usepackage{amsmath, amsthm, amscd, amsfonts}
\usepackage{amssymb}
\usepackage{blkarray}
 \newtheorem{thm}{Theorem}[section]
 \newtheorem{cor}[thm]{Corollary}
 \newtheorem{Lemma}[thm]{Lemma}
 \newtheorem{Prop}[thm]{Proposition}
\theoremstyle{definition}
 \newtheorem{defn}[thm]{Definition}
 \theoremstyle{remark}
  \newtheorem{pro}[thm]{Problem}
 \newtheorem*{ex}{Example}
 \newtheorem{rem}[thm]{Remark}

 \numberwithin{equation}{section}

\usepackage{cite}
\usepackage[colorlinks,linkcolor=blue,citecolor=blue]{hyperref}
\usepackage[numbers]{natbib}
\begin{document}
\date{}
\title[Complex symmetric weighted composition operators]{Complex symmetric weighted composition operators on Dirichlet spaces and Hardy spaces in the unit ball}
\author{\bf Xiao-He Hu, Zi-Cong Yang and Ze-Hua Zhou$^*$}
\address{\newline Xiao-He Hu\newline School of Mathematics, Tianjin University, Tianjin 300354, P.R. China.}
\email{huxiaohe94@163.com}
\address{\newline Zi-Cong Yang\newline School of Mathematics, Tianjin University, Tianjin 300354, P.R. China.}
\email{zicongyang@126.com}
\address{\newline Ze-Hua Zhou\newline School of Mathematics, Tianjin University, Tianjin 300354, P.R. China.}
\email{zehuazhoumath@aliyun.com;zhzhou@tju.edu.cn}
\begin{abstract}
In this paper, we investigate when weighted composition operators acting on Dirichlet spaces $\mathcal{D}(\mathbb{B}_{N})$  are complex symmetric with respect to some special conjugations, and provide some characterizations of Hermitian weighted composition operators on $\mathcal{D}(\mathbb{B}_{N})$. Furthermore, we give a sufficient and necessary condition for $J$-symmetric weighted composition operators on Hardy spaces $H^2(\mathbb{B}_{N})$ to be unitary or Hermitian, then some new examples of complex symmetric weighted composition operators on  $H^2(\mathbb{B}_{N})$ are obtained. We also discuss the normality of complex symmetric weighted composition operators on $H^2(\mathbb{B}_{N})$.
\end{abstract}

\subjclass[2010]{MSC:  47B33, 47B15, 47B38, 32A35, 32A37}
\keywords{weighted composition operator, complex symmetric, Hermitian operator, Dirichlet space, Hardy space}
\thanks{\noindent$^{*}$Corresponding author.
\\
The work was supported in part by the National Natural
Science Foundation of China (Grant Nos. 11771323; 11371276; 11301373).}

\maketitle
\section{Introduction}
Let $\mathbb{B}_{N}$ be the unit ball in $\mathbb{C}_{N}$ and $S_N$ denote the unit sphere. Let $H(\mathbb{B}_{N})$ be
the space of all holomorphic functions on $\mathbb{B}_{N}$. Denote by $S(\mathbb{B}_{N})$ the set of all holomorphic self-maps of $\mathbb{B}_{N}$. The Dirichlet space $\mathcal{D}(\mathbb{B}_{N})$ is defined as
$$
\mathcal{D}(\mathbb{B}_{N})=\{f\in H(\mathbb{B}_{N}):|f(0)|^2+\int_{\mathbb{B}_{N}}|\nabla f(z)|^2d\nu_N(z) < \infty \},
$$
where $d \nu_N(z)$ is the normalized volume measure on $\mathbb{B}_{N}$. The Hardy space $H^2(\mathbb{B}_{N})$  is defined as
$$
H^2(\mathbb{B}_{N})=\{f\in H(\mathbb{B}_{N}):\sup_{0 < r < 1}\int_{S_{N}}|f(r\zeta)|^2d\sigma(\zeta) < \infty \},
$$
where $d\sigma(\zeta)$ is the normalized surface measure on $S_N$.
Recall that for any analytic self-map $\varphi\in S(\mathbb{B}_{N})$ and any analytic function $\psi\in H(\mathbb{B}_{N})$, the weighted composition operator is given by
$$W_{\psi,\varphi}f=\psi\cdot f\circ\varphi.$$
If $\psi\equiv1$, we get the composition operator $C_\varphi$.

In the past five decades, the study of (weighted) composition operator attracted attention of the researchers. It is very interesting to explore how the function theoretic behavior of $\varphi$ affects the properties of $C_\varphi$ on various holomorphic function spaces. For general information about composition operator, we refer the readers to book \cite{CM1} for more details.

An anti-linear operators $C$ on a complex Hibert space $\mathcal{H}$ is called a conjugation if it satisfies the following conditions:
\begin{itemize}
\item[(i)]  involution, i.e. $C^2=I$.
\item[(ii)] isometric, i.e. $ \langle Cx, Cy \rangle = \langle y, x \rangle$ for all $x,y\in\mathcal{H}$.
\end{itemize}
A bounded linear operator $T$ on $\mathcal{H}$ is called complex symmetric if there exists a conjugation $C$ such that $TC=CT^*$ $(CTC=T^*)$. We also say that $T$ is a $C$-symmetric operator.

It is well known that the general study of complex symmetric weighted composition operators on $\mathcal{H}$ are derived from the work of Gracia and Putinar in \cite{GP}, the authors show that every normal operator is complex symmetric. Since then, many significant results about the complex symmetric (weighted) composition operators are obtained. In \cite{GW}, Garica and Wogen got that if an operator $T$ is algebraic of order $2$, then $T$ is complex symmetric, so $C_{\varphi_a}$ is complex symmetric when $\varphi_a$ is an
involutive automorphism. Furthermore, an explicit conjugation  operator $J_a=JW_a$ on $H^2(\mathbb{B}_{N})$  such that $C_{\varphi_a}=J_aC_{\varphi_a}^*J_a$ was given in \cite{N2}. Very recently, in \cite{JK} Jung et al. studied  which combinations of weights $\psi$ and maps of the disk $\varphi$ give rise to complex symmetric weighted composition operators with  respect to classical conjugation $Jf(z)=\overline{f(\overline{z})}$.  In \cite{GZ1}, Gao and Zhou  gave a complete description of complex symmetric composition operators on $H^2(\mathbb{D})$  whose symbols are linear fractional.

Since all conjugation can be considered as a product of a $J$-symmetric unitary operator $U$ and the conjugation $J$. Fatehi in \cite{F} find all unitary weighted composition operators which are $J$-symmetric and consider complex symmetric weighted composition operators with special conjugation $W_{k_a,\varphi_a}J$ on $H^2(\mathbb{D})$. Moreover, a criterion for complex symmetric structure of $W_{\psi,\varphi}$ on $H_\gamma(\mathbb{D})$ (with reproducing kernels $K_w^{\gamma}=(1-\overline{w}z)^{-\gamma}$, where $\gamma\in \mathbb{N}$) was discovered in \cite{LK}. In \cite{ZY}, Yuan and Zhou characterized the adjoint of linear fractional composition operators $C_\varphi$ acting on $\mathcal{D}(\mathbb{B}_{N})$. In \cite{LNN}, the authors showed that no nontrivial normal weighted composition operator exist on the Dirichlet space in the unit disk when $\varphi$ is linear-fractional with fixed point $p\in\mathbb{D}$.

 Motivated by these researches, we attempt to generalize some discussion into more spaces, such as $\mathcal{D}(\mathbb{B}_{N})$,
$H^2(\mathbb{B}_{N})$. The rest of the paper is organized as follows: First, we recall some fundamental definitions and theorems concerning our results. Then we examine the question ``Which weighted composition operator on $\mathcal{D}(\mathbb{B}_{N})$ are complex symmetric?". We prove that $W_{\psi,\varphi}$ is $J$-symmetric ($JC_{Uz}$-symmetric) if and only if $W_{\psi,\varphi}$ is a multiple of the corresponding complex symmetric composition operator $C_\varphi$. In addition, we also show that $C_\varphi$ is  $J$-symmetric if and only if $\varphi(z) = \varphi'(0)z$ for $z\in\mathbb{B}_{N}$ and $\varphi'(0)$ is a symmetric matrix with $||\varphi'(0)||\leq1$. We then provide characterizations of Hermitian  weighted composition operators on $\mathcal{D}(\mathbb{B}_{N})$. Moreover,  we study  when the class of $J$-symmetric weighted composition operator to be  unitary or Hermitian. By providing some sufficient conditions for weighted composition operators to be both unitary and $J$-symmetric, then we get some new examples of complex symmetric weighted composition operators on  $H^2(\mathbb{B}_{N})$. Finally, we discuss the normality of complex symmetric weighted composition operators on $H^2(\mathbb{B}_{N})$.
\section{Preliminaries}

\subsection{Linear fractional map}
\begin{defn} A linear fractional map $\varphi$ of $\mathbb{C}^N$ is a map of the form
$$\varphi(z)=\frac{Az+B}{\langle z,C\rangle+D},$$
where $A = (a_{j,k})$ be an $N \times N$-matrix, $B=(b_j)$, $C =(c_i)$ be $N$-column vectors, $D$ be a complex number, and $\langle . , .\rangle$ indicates the usual Euclidean inner product in  $\mathbb{C}^N$. If $\varphi(\mathbb{B}_{N}) \in \mathbb{B}_{N}$, $\varphi$ is said to be a linear fractional self-map of $\mathbb{B}_{N}$ and signed as $\varphi(z)\in \textrm{LFT} (\mathbb{B}_{N})$.

In this paper, we identify $N \times N$-matrices with linear transformations of $\mathbb{C}_{N}$  via the standard basis of $\mathbb{C}_{N}$.
\end{defn}

\begin{defn}
If $\varphi(z)=\frac{Az+B}{\langle z, C \rangle + D}$ is a linear fractional map, the matrix
$$
 m_{\varphi}=\left(
  \begin{array}{ccc}
  A & B\\
    C^*& D\\
  \end{array}
\right)
$$
will be called a matrix associated with $\varphi$.  If $\varphi(z)\in$ LFT$(\mathbb{B}_{N})$, the adjoint map $\sigma=\sigma_\varphi$ is defined by
$$\sigma(z)=\frac{A^*z-C}{\langle z, -B \rangle + D^*}$$
and the associated matrix of $\sigma$ is
$$
 m_{\sigma}=\left(
  \begin{array}{ccc}
  A^* & -C\\
    -B^*& D^*\\
  \end{array}
\right),
$$
where $A^* = (\overline{a_{ji}})$ denote the conjugate transpose matrix of $A$.
\end{defn}
\begin{thm}{\rm ([\citealp{CM2}, Theorem 4])}\label{Th2.7}
If the matrix
$$
 m_{\varphi}=\left(
  \begin{array}{ccc}
  A & B\\
    C^*& D\\
  \end{array}
\right)
$$
is a multiple of an isometry on the Kre$\breve{{\i}}$n space with
$$
J=\left(
  \begin{array}{ccc}
  I & 0\\
   0& -1\\
  \end{array}
\right)
,$$
then $\varphi(z)=\frac{Az + B}{\langle z, C \rangle + D}$ maps the unit ball $\mathbb{B}_N$ onto itself. Conversely, if $\varphi(z)=\frac{Az+B}{\langle z, C \rangle + D}$ is a linear fractional map of the unit ball onto itself, then $m_\varphi$ is a multiple of an isometry.
\end{thm}

Fix a vector $a\in\mathbb{B}_N$, we denote by $\varphi_a: \mathbb{B}_N\rightarrow\mathbb{B}_N$ the linear fractional map
$$\varphi_a(z)=\frac{a-P_az-s_aQ_az}{1-\langle z, a \rangle},$$
where $s_a=\sqrt{1-|a|^2}$, $P_a$ be the orthogonal projection of $\mathbb{C}^n$ onto the complex line generated by $a$ and $Q_a=I-P_a$, or equivalently,
$$\varphi_a(z)=\frac{a-Tz}{1-\langle z, a \rangle},$$
where $T$ is a self-adjoint map depending on $a$. We use Aut$(\mathbb{B}_N)$ to denote the set of all automorphism of $\mathbb{B}_N$.

\subsection{Spaces and weighted composition operator}
In the Dirichlet space $\mathcal{D}(\mathbb{B}_N)$, evaluation at $w$ in the unit ball is given by $f(w)=\langle f, K_w \rangle$ where
$$K_w(z)=1+\ln\frac{1}{1 - \langle z, w \rangle} \ \ \  \textrm{and} \ \ \  ||K_w||^2= 1 + \ln \frac{1}{{1-|w|^2}}.$$
Let $k_w$ be the normalization of $K_w$, then
$k_w(z)=\frac{(1-\ln(1-|w|^2))^{-\frac{1}{2}}}{1-\ln{1-\langle z, w \rangle}}$.
And in $H^2(\mathbb{B}_{N})$ the kernel for evaluation at $w$ is given by
$$K_w(z)=\frac{1}{(1-\langle z, w \rangle)^N} \ \ \ \textrm{and} \ \ \  ||K_w||^2=(1-|w|^2)^{-N},$$
then $k_w=\frac{(1-|w|^2)^{\frac{N}{2}}}{(1-\langle z, w \rangle)^N}$.

Next we list some fundamental  properties of bounded weighted
composition operators.
\begin{Prop}{\rm ([\citealp{GZ2}, Remark 2.4])}\label{prop2.4}
If $\varphi$ is an automorphism of $\mathbb{B}_N$ and $\psi\in A(\mathbb{B}_N)$ where $A(\mathbb{B}_N)$ denotes the set of functions that holomorphic on $\mathbb{B}_N$ and continuous up to the boundary $S_N$, then
$$W_{\psi,\varphi}=M_\psi C_\varphi.$$
\end{Prop}

\begin{Prop}{\rm ([\citealp{GZ2}, Proposition 2.5])}\label{prop2.5}
 Suppose that $W_{\psi,\varphi}$ is bounded on $\mathcal{D}(\mathbb{B}_N)(H^2(\mathbb{B}_N))$, then we have
$$
W_{\psi,\varphi}^*K_w(z)=\overline{\psi(w)}K_{\varphi(w)}(z)
$$
for all $z,w\in \mathbb{B}_N$. In particular, since $C_\varphi=W_{1,\varphi}$, we get $C_\varphi^*K_w(z)=K_{\varphi(w)}(z)$.
\end{Prop}

\begin{thm}{\rm ([\citealp{CM2}, Theorem 16])}\label{Th2.6}
Suppose  $\varphi(z)=\frac{Az+B}{\langle z, C \rangle + D}$ is a linear fractional map of $\mathbb{B}_N$ into itself for which $C_\varphi$ is a bounded operator on $\mathcal{H}$. Let $\sigma(z)=\frac{A^*z-C}{\langle z,-B \rangle + D^*}$ be the adjoint mapping.
Then $C_\sigma $ is a bounded operator on $\mathcal{H}$, $g(z)=(\langle z, -B\rangle + D^*)^{-r}$ and
$h(z)= (\langle z, C \rangle + D)^r$ are in $H^\infty(\mathbb{B}_N)$, and
$$
C_\varphi^*=T_gC_\sigma T_h^*.
$$
\end{thm}
\subsection{Some others notations}
In this section, we first recall the first partial  derivative reproducing kernel on $\mathcal{D}(\mathbb{B}_N)$.
Let $K_a^{D_1}, K_a^{D_2},\ldots, K_a^{D_n}$ denote the kernels for the first partial derivatives at $a$, that is
$$
\langle f,K_a^{D_j}\rangle=\frac{\partial f}{\partial z_j}(a), \ \ \ \ \ \ \ \   j=1,2,\ldots,n.
$$

It can be shown that
\begin{align}\label{Eq2.1}
K_a^{D_j } = \frac{z_j} {1-\langle z,a \rangle}, \ \ \ \ \ \ \ \ \  j=1,2,\ldots,n.
\end{align}

Continue to the kernels for the first partial derivatives at $a$, we have
\begin{align*}
\langle f , W_{\psi,\varphi}^*K_a^{D_k} \rangle &=\langle W_{\psi,\varphi}f,K_a^{D_k}\rangle\\
&=\frac{\partial\psi(a)}{\partial z_k}f(\varphi(a))+\psi(a)\sum_{j=1}^{n}\frac{\partial f(\varphi(a))}{\partial \varphi_j(z)}\frac{\partial\varphi_j}{\partial z_k}(a)\\
&=\langle f,\overline{\frac{\partial\psi(a)}{\partial z_k}}K_{\varphi(a)}+\overline{\psi(a)}\sum_{j=1}^{n} \overline{\frac{\partial\varphi_j(a)}{\partial z_k}}K_{\varphi(a)}^{D_j} \rangle
\end{align*}
for any $k=1,2,\ldots,n$ and $f\in \mathcal{D}(\mathbb{B}_{N}).$

Thus, we have
\begin{align}\label{Eq2.2}
W_{\psi,\varphi}^*K_a^{D_k} = \overline{\frac{\partial\psi(a)}{\partial z_k}}K_{\varphi(a)}+\overline{\psi(a)}\sum_{j=1}^{n} \overline{\frac{\partial\varphi_j(a)}{\partial z_k}}K_{\varphi(a)}^{D_j}
\end{align}
for any $k=1,2,\ldots,n$ and $f\in \mathcal{D}(\mathbb{B}_{N}).$

In the same way,  we will write
$\langle f, K_a^{D_{i,j}} \rangle = \frac{\partial^2 f}{\partial z_i z_j}$ and
\begin{align}\label{Eq2.3}
 K_a^{D_{i,j}}= \frac{z_i z_j} {(1-\langle z,a \rangle)^2}
\end{align}
for $i,j=1,2,\ldots,n.$
Also, we find that
\begin{align}\nonumber
\langle f , W_{\psi,\varphi}^* K_a^{D_{11}} \rangle \nonumber
&= \langle W_{\psi,\varphi}f , K_a^{D_{11}} \rangle\\ \nonumber
&=\frac{\partial(\frac{\partial \psi(z)} {\partial z_1} f(\varphi(z)) + \psi(z)\sum\limits_{j=1}^{n}\frac{\partial f(\varphi(z))}{\partial z_j}\frac{\partial \varphi_j}{\partial z_1})} {\partial z_1}\Big{|}_{z=a}\\ \nonumber
&=\{\frac{\partial^2 \psi(z)}{\partial z_1 ^2} f(\varphi(z)) + 2 \frac{\partial \psi(z)}{\partial z_1}\sum\limits_{j=1}^{n}\frac{\partial f(\varphi(z))} {\partial z_j} \frac{\partial \varphi_j}{\partial z_1} \\ \nonumber
&+ \psi(z) \frac{\partial (\sum\limits_{j=1}^{n}\frac{\partial f(\varphi(z))} {\partial z_j} \frac{\partial \varphi_j}{\partial z_1}) }{\partial z_1}\}\Big{|}_{z=a}\\ \nonumber
&=\frac{\partial^2 \psi(a)}{\partial z_1 ^2} f(\varphi(a)) + 2 \frac{\partial \psi(a)}{\partial z_1}\sum\limits_{j=1}^{n}\frac{\partial f(\varphi(a))} {\partial z_j} \frac{\partial \varphi_j(a)}{\partial z_1}\\ \nonumber
&+ \psi(a) \sum\limits_{j=1}^{n} \frac{\partial^2 \varphi_j (a)}{ \partial z_1^2} \frac{ \partial f(\varphi(a))}{ \partial z_j}
+ \psi(a)\frac{ \partial \varphi_1(a)}{ \partial z_1} (\sum\limits_{j=1}^{n} \frac{\partial^2 f(\varphi(a))}{  \partial z_1 \partial z_j}\frac{\partial \varphi_j(a)}{\partial z_1}) \\ \nonumber
& +\ldots + \psi(a) \frac { \partial \varphi_n(a)} { \partial z_1} (\sum\limits_{j=1}^{n} \frac{\partial^2 f(\varphi(a))}{ \partial z_n \partial z_j} \frac{\partial \varphi_j(a)}{\partial z_1})\\ \nonumber
&=\langle f, \overline{\frac{\partial^2 \psi(a)}{\partial z_1 ^2}} K_{\varphi(a)}  +  2\overline{\frac{\partial \psi(a)}{\partial z_1}}\sum\limits_{j=1}^{n}\overline{ \frac{\partial \varphi_j(a)}{\partial z_1}} K_{\varphi(a)}^{D_j}\\ \nonumber
&+ \overline{\psi(a)} \sum\limits_{j=1}^{n} \overline{\frac{\partial^2 \varphi_j (a)}{ \partial z_1^2}}  K_{\varphi(a)}^{D_j}
+ \overline{\psi(a)} \overline{ \frac{ \partial \varphi_1(a)} { \partial z_1} }(\sum\limits_{j=1}^{n} \overline{\frac{\partial \varphi_j(a)}{\partial z_1}} K_{\varphi(a)}^{D_{1j}})\\ \nonumber
& +  \ldots  +   \overline{\psi(a)} \overline{ \frac{ \partial \varphi_n(a)} { \partial z_1} }(\sum\limits_{j=1}^{n} \overline{\frac{\partial \varphi_j(a)}{\partial z_1}} K_{\varphi(a)}^{D_{nj}}) \rangle.
\end{align}

Therefore,
\begin{align}
W_{\psi,\varphi}^* K_a^{D_{11}} \nonumber
&=\overline{\frac{\partial^2 \psi(a)}{\partial z_1 ^2}} K_{\varphi(a)}  +  2\overline{\frac{\partial \psi(a)}{\partial z_1}}\sum\limits_{j=1}^{n}\overline{ \frac{\partial \varphi_j(a)}{\partial z_1}} K_{\varphi(a)}^{D_j}\\ \nonumber
& + \overline{\psi(a)} \sum\limits_{j=1}^{n} \overline{\frac{\partial^2 \varphi_j (a)}{ \partial z_1^2}}  K_{\varphi(a)}^{D_j}
+ \overline{\psi(a)} \overline{ \frac{ \partial \varphi_1(a)} { \partial z_1} }(\sum\limits_{j=1}^{n} \overline{\frac{\partial \varphi_j(a)}{\partial z_1}} K_{\varphi(a)}^{D_{1j}}) \\\label{Eq2.4}
& + \ldots   +  \overline{\psi(a)} \overline{ \frac{ \partial \varphi_n(a)} { \partial z_1} }(\sum\limits_{j=1}^{n} \overline{\frac{\partial \varphi_j(a)}{\partial z_1}} K_{\varphi(a)}^{D_{nj}}).
\end{align}
\section{Complex symmetric composition operators on $\mathcal{D}(\mathbb{B}_N)$}

\subsection{Complex symmetric composition operators on $\mathcal{D}$}
Let us start by characterize composition operators on the Dirichlet space in the unit disk which are complex symmetric with respect to the conjugation $Jf(z)=\overline{f(\bar{z})}$. First  we  give the following theorem, which limits the kinds
of maps that can induce complex symmetric composition operators on $\mathcal{D}(\mathbb{D})$.
\begin{thm}
Let $\varphi$ be an analytic self-map of the unit disk. Suppose that $C_\varphi$ is a complex symmetric operator on $\mathcal{D}(\mathbb{D})$, then $\varphi$ has a fixed point in the unit disk.
\end{thm}
\textit{Proof.}
Since $C_\varphi$ is complex symmetric with non-empty point spectrum. [\citealp{N1}, Proposition 3.1] shows that $C_\varphi$ is not hypercyclic. By [\citealp{YR}, Theorem 1.1], we obtain $\varphi$ has a fixed point in  the unit disk.
\qed
\begin{thm}
Let $\varphi$ be an analytic self-map of the unit disk. Then $C_\varphi$ is $J$-symmetric on $\mathcal{D}$ if and only if $C_\varphi$ is normal.
\end{thm}
{\textit{Proof.}}
Since  $C_\varphi$ is $J$-symmetric, it follows from [\citealp{GH}, Proposition 2.4] that $\varphi(z)=az$ for some $|a|\leqslant1$, and so $C_\varphi$ is normal.

For the converse, suppose $C_\varphi$  is normal on $\mathcal{D}$, by [\citealp{CM1}, Theorem 8.2], we conclude that $\varphi(z)=az$ with $|a|\leqslant1$. An easy calculation gives $C_\varphi JK_w(z)=JC_\varphi^*K_w(z)$, for all $z,w \in \mathbb{D}$ and hence $C_\varphi$ is $J$-symmetric.
\qed
\begin{rem}
Indeed, every normal operator is complex symmetric, so it is natural to ask ``if there any complex symmetric but not normal composition operator $C_{\varphi}$ whose symbol is not a constant and also not involution?"
\end{rem}

\subsection{Complex symmetric weighted composition operators on $\mathcal{D}(\mathbb{B}_{N})$}
Following this idea, one is interested in determining whether the $J$-symmetric of composition operator is equivalent to its normality for the dimension greater than 1. Now, we begin with the theorem that gives the sufficient and necessary condition for weighted composition operators $W_{\psi,\varphi}$ to be $J$-symmetric.
\begin{thm}\label{Th3.4}
Let $\varphi$ be an analytic self-map of $\mathbb{B}_N$  and $\psi$ be an analytic function on $\mathbb{B}_{N}$ for which $W_{\psi,\varphi}$ is bounded on $\mathcal{D}(\mathbb{B}_{N})$. Then $W_{\psi,\varphi}$ is complex symmetric with conjugation $J$ if and only if $\psi(z)=c$ and $\varphi(z)=\varphi'(0)z,$ where $c$ is a constant and $\varphi'(0)$ is a symmetric matrix with $||\varphi'(0)||\leq1$.
\end{thm}
\textit{proof}.
 If $W_{\psi,\varphi}$ is complex symmetric with conjugation $J$, then we have
$$W_{\psi,\varphi}JK_w(z)=JW_{\psi,\varphi}^*K_w(z)$$
 for all $z,w \in \mathbb{B}_N$, which implies that
 \begin{equation} \label{Eq3.001}
\psi(z)(1+\ln\frac{1}{1-\langle \varphi(z),\overline{w}\rangle})=\psi(w)(1+\ln\frac{1}{1-\langle z, \overline{\varphi(w)} \rangle}).
\end{equation}

 Putting $w=0$ in Equation (\ref{Eq3.001}), then we have
 \begin{equation} \label{Eq3.002}
 \psi(z)=\psi(0)+\psi(0)\ln\frac{1}{1-\langle z, \overline{\varphi(0)}\rangle}.
\end{equation}

Subsituting the formula for $\psi(z)$ into Equation (\ref{Eq3.001}), we obtain
\begin{align}\nonumber
&(1-\ln(1-\langle z, \overline{\varphi(0)}\rangle))(1-\ln(1-\langle \varphi(z), \overline{w}\rangle))\\ \label{Eq3.003}
&=(1-\ln(1-\langle w, \overline{\varphi(0)}\rangle))(1-\ln(1-\langle z, \overline{\varphi(w)}\rangle)).
\end{align}
Taking partial derivate with respect to $w_1$ on the both sides of Equation (\ref{Eq3.003}), we get
 \begin{align} \nonumber
&(1-\ln(1-\langle z, \overline{\varphi(0)}\rangle))\frac{\varphi_1(z)}{1-\langle\varphi(z),\overline{w}\rangle}\\ \nonumber
&=(1-\ln(1-\langle w, \overline{\varphi(0)}\rangle))\frac{\frac{\partial\varphi_1(w)}{\partial w_1}z_1+\ldots+\frac{\partial\varphi_n(w)}{\partial w_1}z_n}{1-\langle z, \overline{\varphi(w)}\rangle}\\ \nonumber
&+(1-\ln(1-\langle z, \overline{\varphi(w)}\rangle))\frac{\varphi_1(0)}{1-\langle w, \overline{\varphi(0)}\rangle}.
\end{align}
Setting $w=0$ in the above equation, we have
\begin{align} \nonumber
\varphi_1(z)=\varphi_1(0) + \frac{\frac{\partial\varphi_1(0)}{\partial w_1}z_1+\ldots+\frac{\partial\varphi_n(0)}{\partial w_1}z_n}{(1-\langle z, \overline{\varphi(0)}\rangle)(1-\ln(1-\langle z, \overline{\varphi(0)}\rangle))}.
\end{align}
Similarly, we get
$$\varphi_k(z)=\varphi_k(0)+\frac{\frac{\partial\varphi_1(0)}{\partial w_k}z_1+\ldots+\frac{\partial\varphi_n(0)}{\partial w_k}z_n}{(1-\langle z, \overline{\varphi(0)}\rangle)(1-\ln(1-\langle z, \overline{\varphi(0)}\rangle))}$$
for $k=1,2,\ldots,n$. So we have
\begin{align}\label{Eq3.004}
\varphi(z)=\varphi(0)+\frac{\varphi'(0)^Tz}{(1-\langle z, \overline{\varphi(0)}\rangle)(1-\ln(1-\langle z, \overline{\varphi(0)}\rangle))},
\end{align}
here $\varphi'(0)^T$ denote the transpose matrix of $\varphi'(0)$.

Next, we claim  that $\varphi(0)=0$ when $W_{\psi,\varphi}$ is complex symmetric with conjugation $J$. Note that
\begin{equation} \label{Eq3.005}
 W_{\psi,\varphi}J \left(
  \begin{array}{ccc}
  K_a^{D_{11}}\\
  \vdots\\
  K_a^{D_{nn}}
\end{array} \right) = JW_{\psi,\varphi}^*\left(
  \begin{array}{ccc}
  K_a^{D_{11}}\\
  \vdots\\
  K_a^{D_{nn}}
\end{array}\right).
\end{equation}

Putting $a=0$ in Equation (\ref{Eq3.005}), and by Equation (\ref{Eq2.3}), we have
\begin{equation} \label{Eq3.006}
 W_{\psi,\varphi} \left(
  \begin{array}{ccc}
  z_1^2\\
  \vdots\\
  z_n^2
\end{array}\right)=JW_{\psi,\varphi}^*\left(
  \begin{array}{ccc}
  K_0^{D_{11}}\\
 \vdots\\
  K_0^{D_{nn}}
\end{array}\right).
\end{equation}
It follows from  Equation (\ref{Eq3.006})  that

\begin{equation} \label{Eq3.007}
\psi(z) \varphi_1^2(z) = J W_{\psi,\varphi}^* K_0^{D_{11}}(z).
\end{equation}

Since we obtain a precise formula for $\psi(z) \ \textrm{and} \  \varphi(z)$ when $W_{\psi,\varphi}$ is complex symmetric with conjugation $J$, thus by Equation (\ref{Eq3.002}) and (\ref{Eq3.004}) we get
\begin{align} \nonumber
&\psi(z) \varphi_1^2(z)\\ \nonumber
&=\psi(0)\varphi_1^2(0) ( 1 + \ln \frac {1} {1-\langle z, \overline{\varphi(0)} \rangle} ) + 2 \psi(0) \varphi_1(0) \sum\limits_{j=1}^{n} \frac{ \partial \varphi_j(0)} {\partial z_1} \frac{z_j} {1-\langle z, \overline{\varphi(0)} \rangle }+\\ \nonumber
& \psi(0) \frac{ \partial \varphi_1(0)} {\partial z_1} \sum\limits_{j=1}^{n} \frac{ \partial \varphi_j(0)} {\partial z_1} \frac{z_1 z_j} {(1-\langle z, \overline{\varphi(0)} \rangle )^2 (1- \ln (1-\langle z, \overline{\varphi(0)} \rangle)) } + \ldots +\\ \label{Eq3.008}
&\psi(0) \frac{ \partial \varphi_n(0)} {\partial z_1} \sum\limits_{j=1}^{n} \frac{ \partial \varphi_j(0)} {\partial z_1} \frac{z_n z_j} {( 1 - \langle z, \overline{\varphi(0)} \rangle )^2 (1- \ln (1-\langle z, \overline{\varphi(0)} \rangle)) }.
\end{align}

On the other hand, by Equation (\ref{Eq2.4}), we have
\begin{align} \nonumber
&JW_{\psi,\varphi}^* K_0^{D_{11} }(z)\\ \nonumber
&=\frac{\partial^2 \psi(0)}{\partial z_1 ^2} K_{\overline{\varphi(0)}}(z)  +  2\frac{\partial \psi(0)}{\partial z_1}\sum\limits_{j=1}^{n} \frac{\partial \varphi_j(0)}{\partial z_1} K_{\overline{\varphi(0)}}^{D_j}(z)  +  \psi(0) \sum\limits_{j=1}^{n} \frac{\partial^2 \varphi_j (0)}{ \partial z_1^2} K_{\overline{\varphi(0)}}^{D_j} (z) + \\ \label{Eq3.009}
&  \psi(0)  \frac{ \partial \varphi_1(0)} { \partial z_1} (\sum\limits_{j=1}^{n} \frac{\partial \varphi_j(0)}{\partial z_1} K_{\overline{\varphi(0)}}^{D_{1j}}(z) ) +
\ldots   +  \psi(0)  \frac{ \partial \varphi_n(0)} { \partial z_1} (\sum\limits_{j=1}^{n} \frac{\partial \varphi_j(0)}{\partial z_1} K_{\overline{\varphi(0)}}^{D_{nj}}(z) ).
\end{align}

By Equation (\ref{Eq3.002}), we have
\begin{align}\nonumber
\frac{\partial \psi(z)}{\partial z_1} = \psi(0) \frac{\varphi_1(0)}{1-\langle z, \overline{\varphi(0)} \rangle}
\end{align}
and
\begin{align}\nonumber
\frac{\partial^2 \psi(z)}{\partial z_1^2} = \psi(0) \frac{\varphi_1^2(0)}{(1-\langle z, \overline{\varphi(0)} \rangle)^2}.
\end{align}
Thus,
\begin{align} \label{Eq3.0010}
\frac{\partial \psi(0)}{\partial z_1} = \psi(0)\varphi_1(0)
\end{align}
and
\begin{align} \label{Eq3.0011}
\frac{\partial^2 \psi(0)}{\partial z_1^2} = \psi(0)\varphi_1(0)^2.
\end{align}

Then by Equation (\ref{Eq3.004}), a calculation gives
\begin{align}\nonumber
&\frac {\partial^2 \varphi_1(z)} {\partial z_1^2} \\ \nonumber
&= \frac{-(\frac{\partial\varphi_1(0)}{\partial z_1}z_1+\ldots+\frac{\partial\varphi_n(0)}{\partial z_1}z_n ) \frac{\partial^2 F(z) }{ \partial z_1^2}F(z)^2} { {F(z)}^4}-\\ \nonumber
&\frac{\{\frac{\partial\varphi_1(0)}{\partial z_1} F(z)
-{(\frac{\partial\varphi_1(0)}{\partial z_1}z_1 + \ldots + \frac{\partial\varphi_n(0)}{\partial z_1}z_n) \frac{ \partial F(z)}{ \partial z_1}\}2F(z)\frac{\partial F(z)}{ \partial z_1}}} { {F(z)}^4 },
\end{align}
where $F(z) = (1-\langle z, \overline{\varphi(0)}\rangle)(1-\ln(1-\langle z, \overline{\varphi(0)}\rangle))$. Since $\frac{\partial F(0)}{ \partial z_1}=0$, so we have
\begin{align}\label{Eq3.0012}
\frac{\partial^2 \varphi_1(0)}{\partial z_1^2} = 0.
\end{align}

Similarly, we have
\begin{align}\label{Eq3.0013}
\frac{\partial^2 \varphi_j(0)}{\partial z_1^2} = 0
\end{align}
for $j=1,2,\ldots, n.$

Putting all our information together and returning to the Equation (\ref{Eq3.009}), we get
\begin{align} \nonumber
&JW_{\psi,\varphi}^* K_0^{D_{11}}(z)\\ \nonumber
&=\psi(0)\varphi_1(0)^2 (1 + \ln \frac {1} {1-\langle z, \overline{\varphi(0)} \rangle}) + 2 \psi(0) \varphi_1(0) \sum\limits_{j=1}^{n} \frac{ \partial \varphi_j(0)} {\partial z_1} \frac{z_j} {1-\langle z, \overline{\varphi(0)} \rangle }\\ \nonumber
& + \psi(0) \frac{ \partial \varphi_1(0)} {\partial z_1} \sum\limits_{j=1}^{n} \frac{ \partial \varphi_j(0)} {\partial z_1} \frac{z_1 z_j} {(1-\langle z, \overline{\varphi(0)} \rangle )^2  } + \ldots +\\ \label{Eq3.0014}
&\psi(0) \frac{ \partial \varphi_n(0)} {\partial z_1} \sum\limits_{j=1}^{n} \frac{ \partial \varphi_j(0)} {\partial z_1} \frac{z_n z_j} {( 1 - \langle z, \overline{\varphi(0)}\rangle )^2 }.
\end{align}

Combining Equation (\ref{Eq3.008}) and Equation (\ref{Eq3.0014}), we have
\begin{align}  \label{Eq3.0015}
\varphi(0)=0.
\end{align}
Finally, from Equation (\ref{Eq3.002}),  Equation (\ref{Eq3.004}) and  Equation (\ref{Eq3.0015}) we easily deduce that
$$ \psi(z) = \psi(0) = c \ \ \ \ \ \textrm{and} \ \ \ \ \  \varphi(z) = \varphi'(0)z $$
where $c$ is constant and $\varphi'(0)$ is a symmetric matrix with $||\varphi'(0)||\leq1$. Indeed, notice that if $C_{\varphi'(0)z} $ is $J$-symmetric, it is easy to show that $\varphi'(0)=\overline{\varphi'(0)}^*$, that is, $\varphi'(0)$ is a symmetric matrix.

The converse is clear.
\qed
\begin{cor}
Let $\varphi$ be an analytic self-map of $\mathbb{B}_N$, if $C_\varphi$ is bounded on  $\mathcal{D}(\mathbb{B}_{N})$, then $C_\varphi$ is $J$-symmetric on  $\mathcal{D}(\mathbb{B}_{N})$ if and only if $\varphi(z)=\varphi'(0)z$
for $z\in\mathbb{B}_{N}$ and $\varphi'(0)$ is a symmetric matrix  with $||\varphi'(0)||\leq1$.
\end{cor}
 Next we consider the unitary composition operator $C_\varphi$  on $\mathcal{D}(\mathbb{B}_{N})$. Then we can use the unitary composition operator to construct another conjugation operator on  $\mathcal{D}(\mathbb{B}_{N})$.

\begin{Lemma}{\rm ([\citealp{ZY}, Theorem 4.1])}
Let $\varphi$ be an analytic self-map of $\mathbb{B}_{N}$. Then  $C_\varphi$  is unitary on $\mathcal{D}(\mathbb{B}_{N})$  if and only if  $\varphi(z)=Uz$ where $U$ is a unitary matrix.
\end{Lemma}
\begin{Prop} If $U$ is a unitary symmetric matrix, then $JC_{Uz}$ is a conjugation.
\end{Prop}

Using similar proof of Theorem \ref{Th3.4}, we also easily prove the following theorem.

\begin{thm}
Let $\varphi$ be an analytic self-map of $\mathbb{B}_N$  and $\psi$ be an analytic function on $\mathbb{B}_{N}$ for which $W_{\psi,\varphi}$ is bounded on $\mathcal{D}(\mathbb{B}_{N})$. If $W_{\psi,\varphi}$ is complex symmetric with conjugation $JC_{Uz}$ if and only if $\psi(z)=c$ and $\varphi(z)=\varphi'(0)\overline{U}z,$
where $c$ is constant, $\varphi'(0)$ is a symmetric matrix with $||\varphi'(0)||\leq1$ and $\varphi'(0) \overline{U}= \overline{U}\varphi'(0)$.
\end{thm}

\textit{proof}.
Since $W_{\psi,\varphi}$ is complex symmetric with conjugation $JC_{Uz}$, then
$$JC_{Uz}W_{\psi,\varphi}JC_{Uz}=(W_{\psi,\varphi})^*$$
which means that

$$JC_{Uz}W_{\psi,\varphi}=(C_{Uz}W_{\psi,\varphi})^*J.$$

It follows from  Theorem \ref{Th3.4} that
$$\psi(Uz)=\psi(0)$$
and
$$\varphi(Uz)=\varphi'(0)z .$$

Therefore, replace $Uz$ by $z$,  we have
$$\psi(z)=c \ \ \ \ \ \textrm{and} \ \ \ \ \ \ \varphi(z)=\varphi'(0)\overline{U}z,$$
where $c$ is constant, $\varphi'(0)$ is a symmetric matrix with $||\varphi'(0)||\leq1$ and $\varphi'(0)\overline{U} = \overline{U}\varphi'(0)$. In fact, notice that if $C_\varphi$ is $JC_{Uz}$-symmetric, it is  easy to check that  that  $\varphi'(0)\overline{U} = \overline{U}\varphi'(0)$.

The converse direction follows readily from a simple calculation, so we omit the proof.
\qed
\subsection{Hermitian weighted composition operators on $\mathcal{D}(\mathbb{B}_{N})$}
In this section, we will find out the functions $\psi$ and $\varphi$ when $W_{\psi,\varphi}$ are bounded Hermitian weighted composition operators. Not surprisingly, we will prove that no nontrivial Hermitian weighted composition operator exist on  $\mathcal{D}(\mathbb{B}_{N})$.
\begin{thm}
Let $\varphi$ be an analytic self-map of $\mathbb{B}_N$  and $\psi$ be an analytic function on $\mathbb{B}_{N}$ for which $W_{\psi,\varphi}$ is bounded on $\mathcal{D}(\mathbb{B}_{N})$. Then $W_{\psi,\varphi}$ is a Hermitian weighted composition operator on $\mathcal{D}(\mathbb{B}_{N})$ if and only if
 $$\psi(z)=c \ \ \ \textrm{and}\ \ \ \varphi(z)=\varphi'(0)z$$
where $c$ is a real constant and $\varphi'(0)$ is a Hermitian matrix.

\end{thm}
\textit{Proof.}
Since $W_{\psi,\varphi}$ is a bounded Hermitian weighted operator on $\mathcal{D}(\mathbb{B}_{N})$, then we have
$$W_{\psi,\varphi}K_w(z)=W_{\psi,\varphi}^*K_w(z)$$
for all $z,w\in\mathbb{B}_{N}$, that is
$$\psi(z)K_w(\varphi(z))=\overline{\psi(w)}K_{\varphi(w)}(z)$$
for all $z,w\in\mathbb{B}_{N}$. Thus
\begin{equation} \label{Eq3.02}
\psi(z)(1+\ln\frac{1}{1-\langle\varphi(z),w\rangle})=\overline{\psi(w)}(1+\ln\frac{1}{1-\langle z,\varphi(w)\rangle}).
\end{equation}

Letting $w=0$ in Equation (\ref{Eq3.02}) gives
\begin{equation} \label{Eq3.03}
 \psi(z)=\overline{\psi(0)}(1+\ln\frac{1}{1-\langle z,\varphi(0)\rangle})
\end{equation}
for all $z\in\mathbb{B}_{N}$. Putting $z=0$, we get $\psi(0)=\overline{\psi(0)}$, i.e. $\psi(0)$ is a real number.

Substituting the formula for $\psi(z)$ into  Equation (\ref{Eq3.02}), we obtain
\begin{align}\nonumber
&(1-\ln(1-\langle z,\varphi(0)\rangle))(1-\ln(1-\langle\varphi(z),w\rangle)\\\nonumber
&=(1-\ln(1-\langle\varphi(0),w\rangle))(1-\ln(1-\langle z,\varphi(w)\rangle)).
\end{align}
Taking partial derivative with respect to $\overline{w_1}$, we obtain
\begin{align} \nonumber
&(1-\ln(1-\langle z,\varphi(0)\rangle))\frac{\varphi_1(z)}{1-\langle\varphi(z),w\rangle}\\ \nonumber
&=(1-\ln(1-\langle\varphi(0),w\rangle))\frac{\frac{\partial\overline{\varphi_1(w)}}{\partial \overline{w_1}}z_1+\ldots+\frac{\partial\overline{\varphi_n(w)}}{\partial {\overline{w_1}}}z_n}{1-\langle z,\varphi(w)\rangle}\\ \nonumber
&+(1-\ln(1-\langle z,\varphi(w)\rangle))\frac{\varphi_1(0)}{1-\langle\varphi(0),w\rangle}.
\end{align}

Putting $w=0$ in the above equation, we get
$$\varphi_1(z)=\varphi_1(0)+\frac{\frac{\partial\overline{\varphi_1(0)}}{\partial \overline{w_1}}z_1+\ldots+\frac{\partial\overline{\varphi_n(0)}}{\partial \overline{w_1}}z_n}{(1-\langle z,\varphi(0)\rangle)(1-\ln(1-\langle z,\varphi(0)\rangle))}.$$

Similarly, we get
$$\varphi_k(z)=\varphi_k(0)+\frac{\frac{\partial\overline{\varphi_1(0)}}{\partial \overline{w_k}}z_1+\ldots+\frac{\partial\overline{\varphi_n(0)}}{\partial \overline{w_k}}z_n}{(1-\langle z,\varphi(0)\rangle)(1-\ln(1-\langle z,\varphi(0)\rangle))}$$
for $k=1,2,\ldots,n$. Therefore we obtain
\begin{equation}\label{Eq3.04}
\varphi(z)=\varphi(0)+\frac{\varphi'(0)^*z}{(1-\langle z,\varphi(0)\rangle)(1-\ln(1-\langle z,\varphi(0)\rangle))}.
\end{equation}
Taking derivative with respect to $z$, and putting $z=0$, we have
$$\varphi'(0)=\varphi'(0)^*,$$
that is, $\varphi'(0)$ is a Hermitian matrix.

Furthermore,  we obtain $\varphi(0)=0$, the proof is similar to that of Theorem \ref{Th3.4}, so we omit the details. It follows from  Equation (\ref{Eq3.03}) and Equation (\ref{Eq3.04}) that
 $$
 \psi(z)=c \ \ \ \textrm{and} \ \ \ \varphi(z)=\varphi'(0)z
 $$
where $c$ is a  real constant and $\varphi'(0)$ is a Hermitian matrix.

Conversely, if $\psi(z)=c$ and $\varphi(z)=\varphi'(0)z,$ where $c$ is a real constant and $\varphi'(0)$ is a Hermitian matrix. It is easy to see that
$$W_{\psi,\varphi}^*K_w(z)=W_{\psi,\varphi}K_w(z)$$
for all $z,w\in\mathbb{B}_{N}$,  which completes the proof of the theorem.
\qed
\begin{cor}
Let $\varphi$ be an analytic self-map of $\mathbb{B}_N$, if $C_\varphi$ is bounded on  $\mathcal{D}(\mathbb{B}_{N})$, then $C_\varphi$ is Hermitian on  $\mathcal{D}(\mathbb{B}_{N})$ if and only if
$$\varphi(z)=\varphi'(0)z$$
for $z\in\mathbb{B}_{N}$ and some Hermitian matrix $\varphi'(0)$ with $||\varphi'(0)||\leq1$.
\end{cor}
\section{Complex symmetric weighted composition operators on $H^2(\mathbb{B}_{N})$}

  In this section, we begin with some results about complex symmetric weighted composition operators with respect to the conjugation $J$. Then we give some  nontrivial sufficient conditions for weighted composition operators to be unitary and $J$-symmetric. Based upon this, we obtain more  new examples of complex symmetric weighted composition operators on $H^2(\mathbb{B}_{N})$. Finally, we characterize the normality of $C_{Uz}J$-symmetric weighted composition operators.
\subsection{ Unitary and Hermitian complex symmetric weighted composition operators on $H^2(\mathbb{B}_{N})$}
 We first point out when the class of $J$-symmetric weighted composition operator is unitary or Hermitian.
\begin{thm}
Let $\psi(z)=\frac{a_1}{(1-\langle z,\overline{a_{0}}\rangle)^{N}}$, and $\varphi(z)=\frac{a_0-Az}{1-\langle z,\overline{a_{0}}\rangle}$, where $a_{1}\neq0, a_{0}\in{\mathbb{B}}_{N}$ and $A$ is a symmetric matrix such that $\varphi$ is a self-map of ${\mathbb{B}}_{N}$. Then $W_{\psi,\varphi}$ is unitary if and only if
$$a_1=\lambda(1-|a_{0}|^2)^{\frac{N}{2}}, \ \overline{A}A-\overline{a_{0}}a^{T}_{0}=(1-|a_{0}|^2)I_{N} \ \textrm{and }\ A\overline{a_{0}}-{a_{0}}=0,$$
where $a_{0}\in \mathbb{B}_{N}-\{0\}$.

In particular, if $a_0=0$, $W_{\psi,\varphi}$ is unitary if and only if $\psi(z)=\lambda$ for some $\lambda\in\mathbb{C}$ with $|\lambda|=1$ and $\varphi(z)=Az$ where $A$ is a  unitary and symmetric matrix.
\end{thm}
\textit{Proof.}
First recall that the adjoint of $\varphi$ is defined by
 $$ \sigma(z)=\frac{\overline{a_0}-A^{*}z}{1-\langle z,a_{0}\rangle}, \ \  z\in \mathbb{B}_{N}.$$

 Our hypotheses show that $W_{\psi,\varphi}$ is $J$-symmetric (see [\citealp{WY}, Theorem 3.1]), then by [\citealp{L}, Corollary 3.6], we have $W_{\psi,\varphi}$ is a unitary operator if and only if $\psi(z)=\lambda k_{\varphi^{-1}(0)}^{N}(z)$ for some complex number $\lambda$ with $|\lambda|=1$ and $\varphi(z)$ is an automorphism of ${\mathbb{B}}_{N}$. Therefore $\varphi^{-1}(0)=\sigma(0)=\overline{a_{0}}=\overline{\varphi(0)}$, $a_1=\lambda(1-|a_{0}|^2)^{\frac{N}{2}}$. Since $\varphi$ is an analytic self-map of $\mathbb{B}_{N}$ and $\varphi$ is one-to-one (also see [\citealp{WY}, Theorem 3.1]), due to the Theorem \ref{Th2.7}, we have $\varphi(z)\in$ Aut$(\mathbb{B}_{N})$ if and only if
$$
 m_{\varphi}=\left(
  \begin{array}{ccc}
  -A & a_0\\
    -\overline{a_0}^*& 1\\
  \end{array}
\right)
$$
is a non-zero multiple of Kre$\breve{{\i}}$n isometry on Kre$\breve{{\i}}$n space $\mathbb{C}^{N+1}$. So we have
 $$
 |k|^2\left(\begin{array}{cc}
 -A^*&-\overline{a_0}\\
a_0^*&1\\
\end{array} \right)
\left( \begin{array}{cc}
I_N& 0\\
0&-1\\
\end{array} \right)
\left( \begin{array}{cc}
-A & a_0\\
    -\overline{a_0}^*& 1\\
\end{array}
\right)=
\left( \begin{array}{cc}
I_N & 0\\
0& -1\\
\end{array}
\right).
$$

Then we can obtain from direct computations that
\begin{equation} \label{4.1}
|k|^2(A^*A-\overline{a}_0 \overline{a}_0^*)=I_N,
\end{equation}
\begin{equation} \label{4.2}
|k|^2(-A^*a_0+\overline{a}_0)=0,
\end{equation}
\begin{equation} \label{4.3}
|k|^2(a_0^*a_0-1)=-1.
\end{equation}

If $a_0=0$ in Equation (\ref{4.3}), we get
$$\psi(z)=a_1=\lambda$$
for some $\lambda\in\mathbb{C}$ with $|\lambda|=1$ and
 $$A^*A=AA^*=I_N.$$
Thus we have
 $\varphi(z)=Az$, where $A$ is a symmetric unitary matrix.

If $a_0\neq0$ in Equation (\ref{4.3}), we get $|k|^2=\frac{1}{1-|a_0|^2}$, then substituting the expression of $|k|^2$ into Equation (\ref{4.1}) and (\ref{4.2}) we have
$$\overline{A}A-\overline{a}_0a_0^{T}=(1-|a_0|^2)I_N, \ \ \    A\overline{a}_0-a_0=0.$$

Combining these two cases, we have our conclusion.
\qed
\begin{thm}
Let $\psi(z)=\frac{a_1}{(1-\langle z,\overline{a_{0}}\rangle)^{N}}$, and $\varphi(z)=\frac{a_0-Az}{1-\langle z,\overline{a_{0}}\rangle}$, where $a_{1}\neq0, \  a_{0}\in{\mathbb{B}}_{N}$ and $A$ is a symmetric matrix such that $\varphi$ is a self map of ${\mathbb{B}}_{N}$. Then $W_{\psi,\varphi}$ is Hermitian if and only if $a_1$ is a real number, $a_0$ is a real vector, and $A$ is a real matrix.
\end{thm}
\textit{Proof.}
 We know from the [\citealp{WY}, Theorem 3.1] that $W_{\psi,\varphi}$ is complex symmetric with conjugation $J$, and note that $W_{\psi,\varphi}$ is Hermitian (i.e. $W_{\psi,\varphi}=W_{\psi,\varphi}^*$) if and only if
 $$W_{\psi,\varphi} JK_w(z)=JW_{\psi,\varphi} K_w(z)$$
 for any $z,w\in \mathbb{B}_N$, which implies
 \begin{equation} \label{Eq4.4}
\frac{a_1}{(1-\langle z,\overline{a_{0}}\rangle)^{N}(1-\langle \varphi(z),\overline{w}\rangle)^N}=\frac{\overline{a_1}}{(1-\langle\overline{a_{0}},\overline{z}\rangle)^{N}(1-\langle w,\varphi(\overline{z})\rangle)^N}.
\end{equation}
Putting $w=0$ in Equation (\ref{Eq4.4}) we have
$$\frac{a_1}{(1-\langle z,\overline{a_{0}}\rangle)^{N}}=\frac{\overline{a_1}}{(1-\langle\overline{a_{0}},\overline{z}\rangle)^{N}}.$$
Thus, $a_1$ is a real number and $a_0$ is a real vector.

Combining these with Equation (\ref{Eq4.4}) we get
$$\frac{1}{(1-\langle \varphi(z),\overline{w}\rangle)^N}=\frac{1}{(1-\langle w,\varphi(\overline{z})\rangle)^N}.$$
Hence, we have $\varphi(z)=\overline{\varphi(\overline {z})}$, that is
$$\frac{a_0-Az}{1-\langle z,a_0\rangle}=\frac{a_0-\overline{A}z}{1-\langle z,a_0\rangle}.$$
It follows that $A$ is a real matrix, so we complete our proof.
\qed
\subsection{ Some new examples of complex symmetric weighted composition operators on $H^2(\mathbb{B}_{N})$}
In this section, we will give some new examples of complex symmetric weighted composition operators on $H^2(\mathbb{B}_{N})$. For this purpose, we first present some sufficient conditions for weighted composition operators to be  both unitary and $J$-symmetric.

\begin{Prop}\label{Prop4.2}
Let $\Phi$ be an analytic self-map of $\mathbb{B}_{N}$ and $\Psi$ ba an analytic function. If $\Psi(z)=\lambda$,\ $\Phi(z)=Uz$ where $|\lambda|=1$, $U$ is a unitary and symmetric matrix or $\Psi(z)=\mu\frac{(1-|a|^2)^{\frac{N}{2}}}{(1-\langle z,a\rangle)^{N}}$,\ $\Phi(z)=U\frac{a-Tz}{1-\langle z,a\rangle}$, where $\mu \in \mathbb{C}$, $|\mu|=1$, $Ua=\overline{a}$ and $U$
is a symmetric unitary matrix, and $T$ is a self-adjoint map depending on $a$. Then the weighted composition operator $W_{\Psi,\Phi}$ is unitary and $J$-symmetric.
\end{Prop}
\textit{Proof.}
Clearly  if $\Psi(z)=\lambda$ and $\Phi(z)=Uz$, where $|\lambda|=1$, $U$ is a unitary and symmetric matrix, then $W_{\Psi,\Phi}$ is unitary and $J$-symmetric.

Suppose that
$$\Psi(z)=\mu\frac{(1-|a|^2)^{\frac{N}{2}}}{(1-\langle z,a\rangle)^{N}},\  \  \ \ \  \Phi(z)=U\frac{a-Tz}{1-\langle z,a\rangle},$$
where $|\mu|=1$, $Ua=\overline{a}$. From [\citealp{L}, Corollary 3.6], we  see that  $W_{\Psi,\Phi}$ is unitary.
Notice that $Ua=\overline{a}$, so we have
\begin{align} \nonumber
UT&=\sqrt{1-|a|^2}U+(1 - \sqrt{1-|a|^2}) \frac{Ua\overline{a}^T}{|a|^2}\\
&=\sqrt{1-|a|^2}U+(1 - \sqrt{1-|a|^2}) \frac{\overline{a}a^*}{|a|^2}\label{Eq4.05}
\end{align}
Since $U$ is a symmetric unitary matrix,  we can then use Equation (\ref{Eq4.05})  to obtain $UT$ is a symmetric matrix. Then [\citealp{WY}, Theorem 3.1] implies that $W_{\Psi,\Phi}$ is a $J$-symmetric weighted composition operator.
\qed
\begin{pro}
Try to give a sufficient and necessary condition for a weighted composition operator to be both unitary and $J$-symmetric. Note that the case for $H^2(\mathbb{D})$ was solved by Fatehi \cite{F}.
\end{pro}
Before beginning any proofs, we present an example that is simple enough that unitary and $J$-symmetric weighted composition operator $W_{\Psi,\Phi}$ can be carried out concretely.
\begin{ex}
\begin{itemize}
\item[(i)]Let $\Phi$ be an analytic self-map of $\mathbb{B}_{N}$ and $\Psi$ ba an analytic function. If  $\Psi(z)=\mu\frac{(1-|a|^2)^{\frac{N}{2}}}{(1-\langle z,a\rangle)^{N}}$,\ $\Phi(z)=\frac{a-Tz}{1-\langle z,a\rangle}$, where $\mu \in \mathbb{C}$, $|\mu|=1$ and $a\in\mathbb{R}^N \cap \mathbb{B}^N $. Then the weighted composition operator $W_{\Psi,\Phi}$ is unitary and $J$-symmetric.
\item[(ii)]Let $\Phi$ be an analytic self-map of $\mathbb{B}_{N}$ and $\Psi$ ba an analytic function. If  $\Psi(z)=\mu\frac{(1-|a|^2)^{\frac{N}{2}}}{(1-\langle z,a\rangle)^{N}}$,\ $\Phi(z)=U\frac{a-Tz}{1-\langle z,a\rangle}$, where $\mu \in \mathbb{C}$, $|\mu|=1$, $ a=(a_j) $ with $ a_j \neq 0$ for $j=1,2,\ldots,n$ and
    \setcounter{MaxMatrixCols}{20}
\setlength{\arraycolsep}{6pt}
\newcommand*\matrixrightlabel[2][1em]{\text{\makebox[0em][l]{\hspace{#1}$\leftarrow(#2)$}}}

$$
U=\left(
\begin{array}{c}
e^{-2iarga_1} \qquad \qquad \qquad \\
\qquad e^{-2iarga_2} \qquad \qquad \\
\qquad \ddots \qquad \\
\qquad \qquad \qquad e^{-2iarga_n}
\end{array}
\right).
$$
Then the weighted composition operator $W_{\Psi,\Phi}$ is unitary and $J$-symmetric.
\end{itemize}

\setcounter{MaxMatrixCols}{20}
\setlength{\arraycolsep}{6pt}
\newcommand*\matrixrightlabel[2][1em]{\text{\makebox[0em][l]{\hspace{#1}$\leftarrow(#2)$}}}
\end{ex}

 Next, we will use the unitary and $J$-symmetric weighted composition operator in constructing some special conjugations.

\begin{cor}
If $\Psi$ and $\Phi$ satisfy the conditions in Proposition \ref{Prop4.2}, then $W_{\Psi,\Phi}J$ is a conjugation.
\end{cor}
\textit{ Proof. }
The conclusion follows directly from [\citealp{GPP}, Lemma 3.2].
\qed

Next we present a sufficient and necessary condition for weighted composition operators to be $W_{\Psi,\Phi}J$-symmetric, where $\Psi$ and $\Phi$ satisfy the conditions in Proposition \ref{Prop4.2}.
Most of them are direct extensions of results in \cite{F}, so we omit the details.
\begin{thm}\label{Th4.5}
Let $\psi(z)=\frac{a_1}{(1-\langle z,\overline{a_{0}}\rangle)^{N}}$, and $\varphi(z)=\frac{a_0-Az}{1-\langle z,\overline{a_{0}}\rangle}$, where $a_{1}\neq0, a_{0}\in{\mathbb{B}}_{N}$ and $A$ is a symmetric matrix such that $\varphi$ is a self map of ${\mathbb{B}}_{N}$.
\begin{itemize}
\item[(i)]   For  $a\neq0$, the weighted composition operator $W_{\widetilde{\psi},\widetilde{\varphi}}$ is complex symmetric with conjugation $W_{\Psi,\Phi}J$ if and only if $\widetilde{\psi}=\Psi\cdot\psi\circ\Phi$, $\widetilde{\varphi}=\varphi\circ\Phi$.
\item[(ii)]  If $U$ is a unitary symmetric matrix, the weighted composition operator $W_{\widetilde{\psi},\widetilde{\varphi}}$ is complex symmetric with conjugation $C_{Uz}J$ if and only if
$\widetilde{\psi}(z)=\psi(Uz)$ and $\widetilde{\varphi}(z)=\varphi(Uz).$
\end{itemize}

\end{thm}
Due to a result of [\citealp{NST}, Theorem 2.10], Narayan, Sievewright, and Thompson give some examples of linear-fractional, not automorphic maps that induce complex symmetric composition operators on $H^2(\mathbb{D})$. The ideals in \cite{NST} may be adapted to prove when $C_{Az+c}$ is $JW_{\psi_{b},\varphi_{b}}$-symmetric on $H^2(\mathbb{B}_{N})$, where $\psi_{b}(z)= k_b^N(z)=\frac{(1-|b|^2)^{\frac{N}{2}}}{(1-\langle z,b\rangle)^N}$
,\ $\varphi_{b}(z)=\frac{b-Tz}{1-\langle z,b\rangle}$ with $b=(I-A)^{-1}c $, and $T = \sqrt{1-|b|^2}I + (1 - \sqrt{1-|b|^2}) \frac{bb^T}{|b|^2}$.
However, the result can fail when $N > 1$, we do have the following restricted version.
\begin{Lemma}\label{lem2.7}
 Let $\gamma >0$ be given. Suppose $\varphi$ is a linear fractional map of $\mathbb{B}_{N}$  into itself, then $C_\varphi$  is bounded on $H_\gamma$.
\end{Lemma}
\begin{Lemma}
 Let $\sigma(z)=Az+c$ is a linear fractional map of $\mathbb{B}_{N}$ into $\mathbb{B}_{N}$. Then $C_\sigma$ is bounded on $H^2(\mathbb{B}_{N})$ and  we have
 $$C_\sigma^*=T_\psi C_\varphi,$$
 where $\psi(z)=\frac{1}{(1+\langle z,-c\rangle)^N}$, $\varphi(z)=\frac{A^*z}{1+\langle z,-c\rangle}.$
\end{Lemma}
\textit{Proof.}
We know from Lemma \ref{lem2.7} that $C_\sigma$ is bounded on $H^2(\mathbb{B}_{N})$, then by the adjoint formula on $H^2(\mathbb{B}_{N})$ given by Theorem \ref{Th2.6}, we will get this result by a simple computation.
\begin{thm}\label{Th4.8}
Let $A=(a_{i,j})$ be $N\times N$-symmetric matrix, and 1 is not a eigenvalue of $A$, $c=(c_i)$ be $N$-column vectors. Let $b=(I-A)^{-1} c\in \mathbb{B}_{N}\cap \mathbb{R}_{N}$, and $Ab=\lambda b$ for some $\lambda\in\mathbb{C}$, let $\sigma(z)=Az+c$, then $C_{\sigma}$ is $JW_{\psi_{b},\varphi_{b}}$-symmetric.
\end{thm}

\textit{Proof.}
 we need to prove that
 \begin{equation} \label{Eq4.5}
 C_\sigma JW_{\psi_{b},\varphi_{b}}K_w(z)=JW_{\psi_{b},\varphi_{b}}C_\sigma^*K_w(z)
\end{equation}
for all $z, w\in\mathbb{B}_{N} $, where $C_{\sigma}^*=M_\psi C_\varphi=\frac{1}{(1+\langle z,-c\rangle)^N}\frac{A^*z}{1+\langle z,-c\rangle}$, $\psi_{b} (z)=k_b^N(z)=\frac{(1-|b|^2)^{\frac{N}{2}}}{(1-\langle z,b\rangle)^N}$, and $\varphi_{b}(z)=\frac{b-Tz}{1-\langle z,b\rangle}$ with $T$ is a self-adjoint operator depending on $b$.

For the left side of of Equation (\ref{Eq4.5}), we have
\begin{align*}
&C_{Az+c}JW_{\psi_{b},\varphi_{b}}K_w(z)\\
&=C_{Az+c}J\frac{(1-|b|^2)^{\frac{N}{2}}}{(1-\langle z,b\rangle)^N}K_w(\frac{b-Tz}{1-\langle z,b\rangle})\\
&=\frac{(1-|b|^2)^{\frac{N}{2}}}{(1-\langle b,\overline{Az+c}\rangle)^N}\overline{K_w(\frac{b-T(\overline{Az+c})}{1-\langle\overline{Az+c},b\rangle})}\\
&=\frac{(1-|b|^2)^{\frac{N}{2}}}{(1-\langle b,\overline{Az+c}\rangle)^N}\frac{1}{(1-\langle w,\frac{b-T(\overline{Az+c})}{1-\langle\overline{Az+c},\overline{b}\rangle}\rangle)^N}\\ &=\frac{(1-|b|^2)^{\frac{N}{2}}}{(1-\langle b,\overline{Az+c}\rangle-\langle w,b-T\overline{Az}-T\overline{c}\rangle)^N}.
\end{align*}

For the right side of of Equation (\ref{Eq4.5}), we have
\begin{align*}
&JW_{\psi_{b},\varphi_{b}}C_{Az+c}^*K_w(z)\\
&=J\frac{(1-|b|^2)^{\frac{N}{2}}}{(1-\langle z,b\rangle)^N}C_\frac{b-Tz}{1-\langle z,b\rangle} \frac{1}{(1+\langle z,-c\rangle)^N} K_w \Big(\frac{A^*z}{1+\langle z,-c\rangle} \Big) \\
&=\frac{(1-|b|^2)^{\frac{N}{2}}}{(1-\langle b,\overline{z}\rangle)^N} \frac{1}{(1+\langle-c,\frac{b-T\overline{z}}{1-\langle\overline{z},b\rangle}\rangle)^N} \frac{1}{(1-\langle w,\frac{A^*b-A^*T\overline{z}}{1-\langle\overline{z},b\rangle+\langle b-T\overline{z},-c\rangle}\rangle)^N} \\
&=\frac{(1-|b|^2)^{\frac{N}{2}}}{(1-\langle b,\overline{z}\rangle+\langle-c,b\rangle+\langle c,T\overline{z}\rangle-\langle Aw,b\rangle+\langle Aw,T\overline{z}\rangle)^N}.
\end{align*}

To see Equation (\ref{Eq4.5}) holds, by the previous calculations, it is enough to show that
\begin{equation} \label{Eq4.6}
-z^TA^Tb - b^Tw + z^TA^TTw + c^TTw = -z^Tb + z^TTc - b^TAw + z^TTAw.
\end{equation}

Since $T = \sqrt{1-|b|^2}I + (1 - \sqrt{1-|b|^2}) \frac{bb^T}{|b|^2}$ and $Ab = \lambda b$ for some $\lambda \in C$, we have $c = kb$ with $k = 1 - \lambda$ and
\begin{align}
Tc &= \sqrt{1-|b|^2}c + \frac{(1 - \sqrt{1-|b|^2})bb^T}{|b|^2}kb \nonumber \\
&= \sqrt{1-|b|^2}kb + (1 - \sqrt{1-|b|^2})kb  \nonumber \\
&= kb = c .\label{Eq4.7}
\end{align}

By now, if we can verify that $A$ and $T$ are commute, i.e.
\begin{equation} \label{Eq4.8}
AT = TA,
\end{equation}
then by Equations (\ref{Eq4.6}), (\ref{Eq4.7}) and (\ref{Eq4.8}), we get Equation (\ref{Eq4.5}) holds, immediately.

Indeed, we just need to prove $Abb^T = bb^TA$ holds. Since
$$
Abb^T = A \frac{c}{k} \Big(\frac{c}{k}\Big)^T = \Big(\frac{c}{k}- c\Big) \frac{c^T}{k} = \frac{cc^T}{k^2} - \frac{cc^T}{k}
$$
and
$$
bb^TA = \frac{c}{k} \Big(\frac{c}{k}\Big)^TA = \frac{c}{k} \Big(\frac{Ac}{k}\Big)^T = \frac{c}{k} \Big(\frac{c}{k} - c\Big)^T = \frac{cc^T}{k^2} - \frac{cc^T}{k},
$$
where we use the fact that $A$ is a symmetric matrix, so we have $Abb^T = bb^TA$, as desired.
\qed

\subsection{Normality of complex symmetric weighted composition operators on $H^2(\mathbb{B}_{N})$}
We turn next to the problem of identifying the when the class of complex symmetric weighted composition operators mentioned in Theorem \ref{Th4.5} (ii) is normal.
\begin{thm}
Let $\psi(z)=\frac{a_1}{(1-\langle Uz,\overline{a_{0}}\rangle)^{N}}$, and $\varphi(z)=\frac{a_0-AUz}{1-\langle Uz,\overline{a_{0}}\rangle}$, where $a_{1}\neq0, a_{0}\in{\mathbb{B}}_{N}$, $A$ is a symmetric matrix and $U$ is a unitary symmetric matrix such that $\varphi$ is a self map of ${\mathbb{B}}_{N}$.  Then $W_{\psi,\varphi}$ is normal if and only if
$$\overline{UA}AU-\overline{Ua_0}a_0^TU=A\overline{A}-a_0{a_0^*}$$
and
$$\overline{UA}a_0 - \overline{Ua_0} = A\overline{a_0} - a_0.$$
\end{thm}
\textit{Proof.}
 Since $\varphi(z)=\frac{a_0-AUz}{1-\langle Uz,\overline{a_{0}}\rangle}$, the associated matrix of $\varphi$ is
$$
 m_{\varphi}=\left(
  \begin{array}{ccc}
  -AU &a_0\\
   {-a_0^T}U& 1\\
  \end{array}
\right)
$$
and the  adjoint map $\sigma=\sigma_\varphi$ is defined by
 $$\sigma(z)=\frac{(-AU)^*z+\overline{Ua_0}}{1+\langle z,-a_{0}\rangle}$$
 with the associated matrix of $\sigma$ is

 $$
 m_{\sigma}=\left(
  \begin{array}{ccc}
  -(AU)^* &\overline{Ua_0}\\
   {-a_0}^*& 1\\
  \end{array}
\right).
$$

 Therefore we have  $|\varphi(0)|=|a_0|=|Ua_0|=|\sigma(0)|$, $ \psi(z)=a_1K_{\overline{Ua_0}}(z)=a_1K_{\sigma(0)}(z)$, \   by  [\citealp{L}, Proposition 4.6], we see that $W_{\psi,\varphi}$ is normal if and only if $\varphi\circ\sigma=\sigma\circ\varphi$. Since the multiple of $m_{\varphi\circ\sigma}$ give the same linear fractional map, we must have $W_{\psi,\varphi}$ is normal if and only if $m_{\varphi\circ\sigma}=k m_{\sigma\circ\varphi}$ for some $k\neq0$.

 A calculation shows that
 \begin{equation}\label{Eq4.10}
\overline{UA}AU-\overline{Ua_0}a_0^TU=k(A\overline{A}-a_0{a_0^*}),
\end{equation}
\begin{equation}\label{Eq4.11}
-\overline{UA}a_0+\overline{Ua_0}=k(-A\overline{a_0}+a_0),
\end{equation}
\begin{equation}\label{Eq4.12}
-\overline{a_0}^Ta_0+1=k(-a_0^T\overline{a_0}+1).
\end{equation}
Because $|a_0|=|\varphi(0)|\neq 1$, by  Equation (\ref{Eq4.12}), we have $k=1$. Then combining this with  Equation (\ref{Eq4.10}) and (\ref{Eq4.11}), we have $W_{\psi,\varphi}$ is normal if and only if  $\overline{UA}AU-\overline{Ua_0}a_0^TU=A\overline{A}-a_0{a_0^*}$  and $\overline{UA}a_0 - \overline{Ua_0} = A\overline{a_0} - a_0$, as desired.
\qed

\end{document}